\documentclass[12pt]{amsart}

\usepackage{amssymb} 
\usepackage{graphics}
\usepackage{epsf} 
\usepackage{amsmath,epsf} 
\usepackage[latin1]{inputenc} 
\usepackage{amsfonts} 
\usepackage[dvips]{epsfig} 
\numberwithin{equation}{section} 
\newtheorem{theo}{Theorem}[section] 
\newtheorem{prop}[theo]{Proposition} 
\newtheorem{conj}[theo]{Conjecture} 
\newtheorem{obs}[theo]{Observation} 
\newtheorem{lemm}[theo]{Lemma} 
\newtheorem{coro}[theo]{Corollary}

\newtheorem{rema}[theo]{Remark}

\def\col{\!\!\uparrow}

\def\D{{\mathcal D}}

\def\charac{\raise2pt\hbox{$\chi$}}

\newcommand{\la}{{\lambda}}

\newcommand{\al}{{\alpha}} 
 
\newcommand{\be}{{\beta}}

\def\Young#1{\vbox{\smallskip\offinterlineskip 
    \halign{&\vbox{##}\kern-\Thickness\cr #1}}} 
 
\newdimen\Squaresize \Squaresize=12pt 
\newdimen\Thickness \Thickness=.3pt 
\newdimen\Correction \Correction=7pt 
 
\def\Vide#1{\hbox{ 
       \vbox to \Squaresize{\vss 
          \hbox to \Squaresize{\hss#1 \hss}\vss} 
    \hskip-\Correction} 
   \kern-\Thickness}

\def\Carre#1{\hbox{\vrule width \Thickness 
   \vbox to \Squaresize{\hrule height \Thickness\vss 
      \hbox to \Squaresize{\hss#1\hss} 
   \vss\hrule height\Thickness} 
   \unskip\vrule width \Thickness} 
   \kern-\Thickness} 
  
\def\Blk{\Vide{}} 
 
\def\Box#1{\Carre{$\scriptstyle#1$}}

\title[Schur-positivity]{Inequalities between Littlewood-Richardson Coefficients} 
\author{F. Bergeron} 
\address[F. Bergeron]{D\'epartement de Math\'ematiques\\ Universit\'e 
  du Qu\'ebec \`a Montr\'eal\\ Montr\'eal, Qu\'ebec, H3C 3P8, CANADA} 
\author{R. Biagioli} 
\author{M. Rosas} 
\address[R. Biagioli, M. Rosas]{LaCIM\\ Universit\'e 
  du Qu\'ebec \`a Montr\'eal\\ Montr\'eal, Qu\'ebec, H3C 3P8, CANADA} 
\email{bergeron.francois@uqam.ca, biagioli@math.uqam.ca, rosas@math.uqam.ca} 
 
\thanks{F. Bergeron is supported in part by NSERC and 
  FQRNT}

\begin{document}
\begin{abstract}
We prove that a conjecture of Fomin, Fulton, Li, and Poon, associated to ordered pairs of partitions,  holds for many infinite families of such pairs. We also show that  the  bounded height case can be reduced to checking that the conjecture holds for a finite number of pairs, for any given height. Moreover, we propose a natural generalization of the conjecture to the case of skew shapes.
\end{abstract}  
 \maketitle
 \parskip=0pt
\tableofcontents
\parskip=8pt  
\section{Introduction}  
 
In the course of their study of Horn type inequalities for eigenvalues and singular values of complex matrices, Fomin, Fulton, Li, and Poon \cite{fulton} come up with a very  
interesting conjecture concerning the Schur-positivity of special differences of products of Schur functions. More precisely, they consider differences of the form
  $$s_{\mu^*}s_{\nu^*}-s_{\mu}s_{\nu},$$
   where $\mu^*$ and $\nu^*$ are partitions constructed from an ordered pair of partitions $\mu$ and $\nu$ through a seemingly strange procedure at first glance. In our presentation, their transformation $(\mu,\nu) \mapsto (\mu^{*},\nu^*)$ on ordered pairs of partitions, will rather be denoted 
\begin{equation}\label{defstar} 
(\mu,\nu) \longmapsto {(\mu,\nu)}^*= (\lambda(\mu,\nu),\rho(\mu,\nu)) 
\end{equation} 
and will be called the $*$-operation. As we shall see, this change of notation is essential in order to simplify the presentation of the many nice combinatorial properties of this operation. On the other hand, it underlines that both entries, $\lambda$ and $\rho$ of the image ${(\mu,\nu)}^*$ of $(\mu,\nu)$, actually depend on both $\mu$ and $\nu$. 
 
With this slight change of notation, the original definition of the $*$-operation is as follows.   Let $\mu=(\mu_1, \mu_2, \ldots,\mu_n)$ and $\nu=(\nu_1, \nu_2, \ldots,\nu_n)$ two partitions with the same number of parts, allowing zero parts.
From these, two new partitions $\lambda(\mu,\nu)=(\lambda_1, \lambda_2, \ldots,\lambda_n)$ and $\rho(\mu,\nu)= (\rho_1, \rho_2, \ldots,\rho_n)$ are constructed as follows 
\begin{equation}\label{definition}
\begin{array}{lcl}
   \lambda_k&:= &\mu_k-k + \# \{j\ |\ 1\leq j\leq n,\  \nu_j-j\geq \mu_k-k\};\\  
   \rho_j&:=&\nu_j-j+1+ \# \{k\ |\ 1\leq k\leq n,\  \mu_k-k>\nu_j-j\}. 
 \end{array}
 \end{equation}

\noindent Although this definition does not make it immediately clear, both $\lambda(\mu,\nu)$ and $\rho(\mu,\nu)$ are truly partitions, and they are such that
\[ |\lambda(\mu,\nu)|\,+\,|\rho(\mu,\nu)|=|\mu|\,+\,|\nu|,\]
where as usual $|\mu|$ denotes the sum of the parts of $\mu$. 

Recall that the product of two Schur  functions can always be expanded as a linear combination
   $$s_\mu s_\nu=\sum_\theta c_{\mu\,\nu}^\theta  s_\theta,$$ 
 of Schur functions indexed by partitions $\theta$ of the integer $n=|\mu|+|\nu|$, since these Schur functions constitute a linear basis of the homogeneous symmetric functions of degree $n$. It is a particularly nice feature of this expansion that the coefficients $c_{\mu\,\nu}^\theta$ are always non-negative integers. They are called the {\em Littlewood-Richardson coefficients}. More generally, we say that a symmetric function is {\em Schur positive} whenever the coefficients in its expansion, in the Schur function basis, are all non-negative integers. For more details on symmetric function theory see Macdonald's classical book \cite{macdonald}, whose notations we will mostly follow.
We can then state the following:
   
\begin{conj}[Fomin-Fulton-Li-Poon]\label{fflp} 
For any pair of partitions $(\mu,\nu)$,  if 
    $${(\mu,\nu)}^*=(\lambda,\rho),$$ 
then the symmetric function  
\begin{equation}\label{conj} 
s_{\lambda}s_{\rho}-s_\mu s_\nu 
\end{equation} 
is Schur-positive. 
\end{conj} 
 
\noindent In other words, this says that $c_{\mu\,\nu}^\theta\leq c_{\lambda\,\rho}^\theta$, for all $\theta$ such that $s_\theta$ appears in the expansion  of $s_\mu s_\nu$.

For an example of one of the simplest case of the $*$-operation, let $\mu=(a)$ and $\nu=(b)$, with $a > b$, be two one-part partitions.  In this case, we get 
   $${((a),(b))}^*=(a-1,b+1),$$
so that Conjecture \ref{fflp} corresponds exactly to an instance of the classical Jacobi-Trudi identity: 
\begin{eqnarray*}
    s_{a-1}s_{b+1} - s_a s_b &=& \det \left( \begin{array}{cc}  
                                        s_{a-1} & s_a \\ 
                                        s_b     & s_{b+1} \end{array} \right)\\
                                      &=&s_{a-1,b+1}.
\end{eqnarray*}
In this paper we give a new recursive combinatorial description of the $*$-operation. This recursive description allows us to prove many instances of Conjecture \ref{fflp} and to show that it reduces to checking a finite number of instances for any fixed $\nu$, if we bound the number of parts of $\mu$.  
Moreover
we show how to naturally generalize the conjecture to  pairs of skew partitions.

\section{Combinatorial properties of the $*$-operation and implications}

We first derive some nice combinatorial properties of the transformation $*$. To help in the presentation of these properties, let us introduce some further notation. For any undefined notation we refer to \cite{macdonald}. We often identify a partition with its (Ferrers) diagram. Diagrams are drawn here using the ``French'' convention of ordering parts in decreasing order from bottom to top.
 
\begin{figure*}[ht]   
\centering  
\begin{picture}(0,0)(0,0)
\put(16,26){$\longrightarrow_2$}
\put(56,11){$=$}
\put(213,13){$\Big\uparrow_2$}
\put(230,11){$=$}
\end{picture}
\includegraphics[width=100mm]{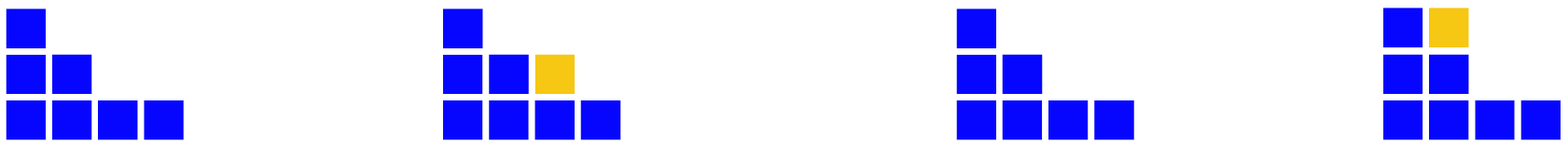}
\caption{ \label{Notation}}   
\end{figure*}

We write $\mu={\overrightarrow{\alpha}}^i$, if the partition $\mu$ is obtained from the partition $\alpha$ by adding one cell in line $i$; and 
$\mu=\alpha\col_k$, if $\mu$ is obtained from $\alpha$ by adding one cell in column $k$. In other words, $\mu={\overrightarrow{\alpha}}^i$ means that $\mu_i=\alpha_i$ for all $i\not=\ell$, and $\mu_\ell=\alpha_\ell+1$. This is illustrated in Figure \ref{Notation} in term of diagrams.
 
 \noindent Observe that, \setlength{\jot}{0pt}
 \begin{eqnarray*}
     \mu={\overrightarrow{\alpha}}^i\qquad &{\rm iff}&\qquad \mu'={\overrightarrow{\alpha'}}^{\mu_i}\label{add_conj_col}\\
     &{\rm iff}&\qquad \mu=\alpha\col_{\mu_i}\label{add_col}\\
                      &{\rm iff}&\qquad \mu'=\alpha'\col_i\label{add_conj}
  \end{eqnarray*}  
\noindent We can now state our recursive description of the $*$-operation.  
 \begin{prop}[Recursive formula] \label{recurrence} For any  partitions $\alpha$ and $\nu$, if $(\lambda,\rho)={(\alpha,\nu)}^*$, then  we have 
      \begin{equation}\label{gauche}
          {({\overrightarrow{\alpha}}^i,\nu)}^*=\begin{cases}
                    (\lambda,{\overrightarrow{\rho}}^j)  & \text{if $j$ is such that  $\nu_j-j=\alpha_i-i$, if any}, \\ \\
                    ({\overrightarrow{\lambda}}^i,\rho) & \text{otherwise}.
           \end{cases} 
       \end{equation} 
Moreover, when in the first case, we have ${\overrightarrow{\rho}}^j=\rho\col_{\mu_i}$.
In a similar manner, for given $\mu$ and $\beta$, if $(\lambda,\rho)={(\mu,\beta)}^*$, then 
      \begin{equation}\label{droite}
          {(\mu,{\overrightarrow{\beta}}^i)}^*=\begin{cases}
                    ({\overrightarrow{\lambda}}^j,\rho)  & \text{if $j$ is such that  $\mu_j-j=\nu_i-i$, if any}, \\ \\
                    (\lambda,{\overrightarrow{\rho}}^i) & \text{otherwise},
           \end{cases} 
       \end{equation} 
  and, when in the first case, we have ${\overrightarrow{\lambda}}^j=\lambda\col_{\nu_i}$.
 \end{prop} 
 
\noindent We can clearly use Proposition \ref{recurrence} to recursively compute $\lambda(\mu,\nu)$ and $\rho(\mu,\nu)$. This is discussed more extensively  in Section \ref{ext_tableau}.  
The actual computation of the $*$-operation can be simplified in view of the following property (see Lemma \ref{conjugate}).  
For any pair of partitions $(\mu, \nu)$, we have
\begin{align}\label{transpose} 
{(\mu,\nu)}^*=(\lambda,\rho)\qquad{\rm iff}\qquad {(\nu',\mu')}^*=(\lambda',\rho'), 
\end{align} 
where, as usual, $\mu^{\prime}$ stands from the conjugate of $\mu$. Using the fact that the involution $\omega$  (which is the linear operator that maps $s_\mu$ to $s_{\mu'}$) is multiplicative, 
it easily follows  that  
\begin{prop}\label{cortranspose}   
Conjecture \ref{fflp} holds for the pair $(\mu,\nu)$ if and only if it holds for the pair $(\nu',\mu')$. 
\end{prop} 
 
In practice, there are many ways to describe the $*$-operation recursively, since we can freely choose how to make partitions grow. It is sometimes convenient to start from the pair $(0, \nu)$,  with $0$ standing for the empty partition, whose image under the  $*$-operation has a simple description.
 
\begin{figure*}[ht] 
\centering  
\begin{picture}(0,0)(0,0)
\put(35,-10){$\nu$}
\put(143,-10){$\overline{\nu}$}
\put(215,-10){$\underline{\nu}$}
\put(96,34){$*$}
\put(90,30){$\longrightarrow$}
\end{picture}
\includegraphics[width=90mm]{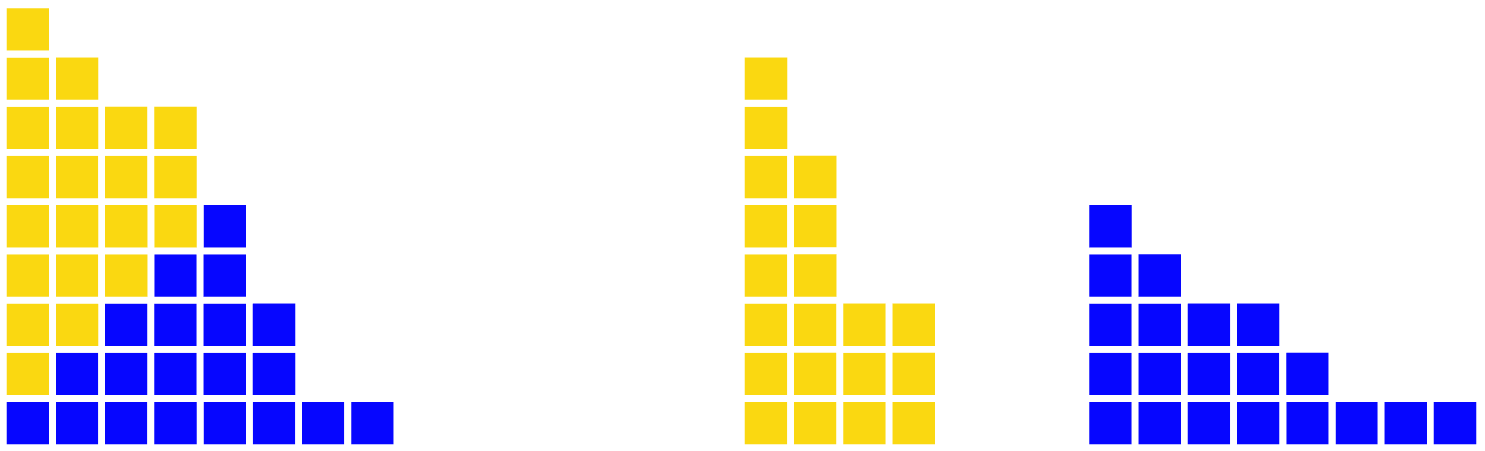}  
\caption{  
\label{original}  } 
\end{figure*}

\begin{lemm}\label{stair} Let $\nu$ be any partition. Then 
\begin{align*} 
\rho(0,\nu)&= (\nu_1, \nu_2-1, \cdots, \nu_k-(k-1))\\ 
\lambda'(0,\nu) &= (\nu'_1-1, \nu'_2-2, \cdots, \nu'_k-k)
\end{align*} 
where $k=\max \{ i : \nu_i-(i-1) \ge 1\}$,.
\end{lemm} 
\noindent We will sometimes denote respectively $\overline{\nu}$ and $\underline{\nu}$ the partitions $\lambda(0,\nu)$ and $\rho(0,\nu)$.
For example if $\nu=866554421$, then 
    $$\overline{\nu}=44432211\qquad {\rm and} \qquad \underline{\nu}=85421$$ 
 as is  illustrated in Figure \ref{original}. In Section \ref{ext_tableau} we elaborate on the various ways that Proposition \ref{recurrence} 
can be used to compute the $*$-operation. This gives rise to a $*$-operation on pairs of Young tableaux.
 
\begin{figure*}[ht] 
\centering 
\begin{picture}(0,0)(0,0)
\put(141,289){$*$}
\put(141,204){$*$}
\put(141,114){$*$}
\put(141,24){$*$}
\put(135,285){$\longrightarrow$}
\put(135,200){$\longrightarrow$}
\put(135,110){$\longrightarrow$}
\put(135,20){$\longrightarrow$}
\end{picture}  
\includegraphics[width=100mm]{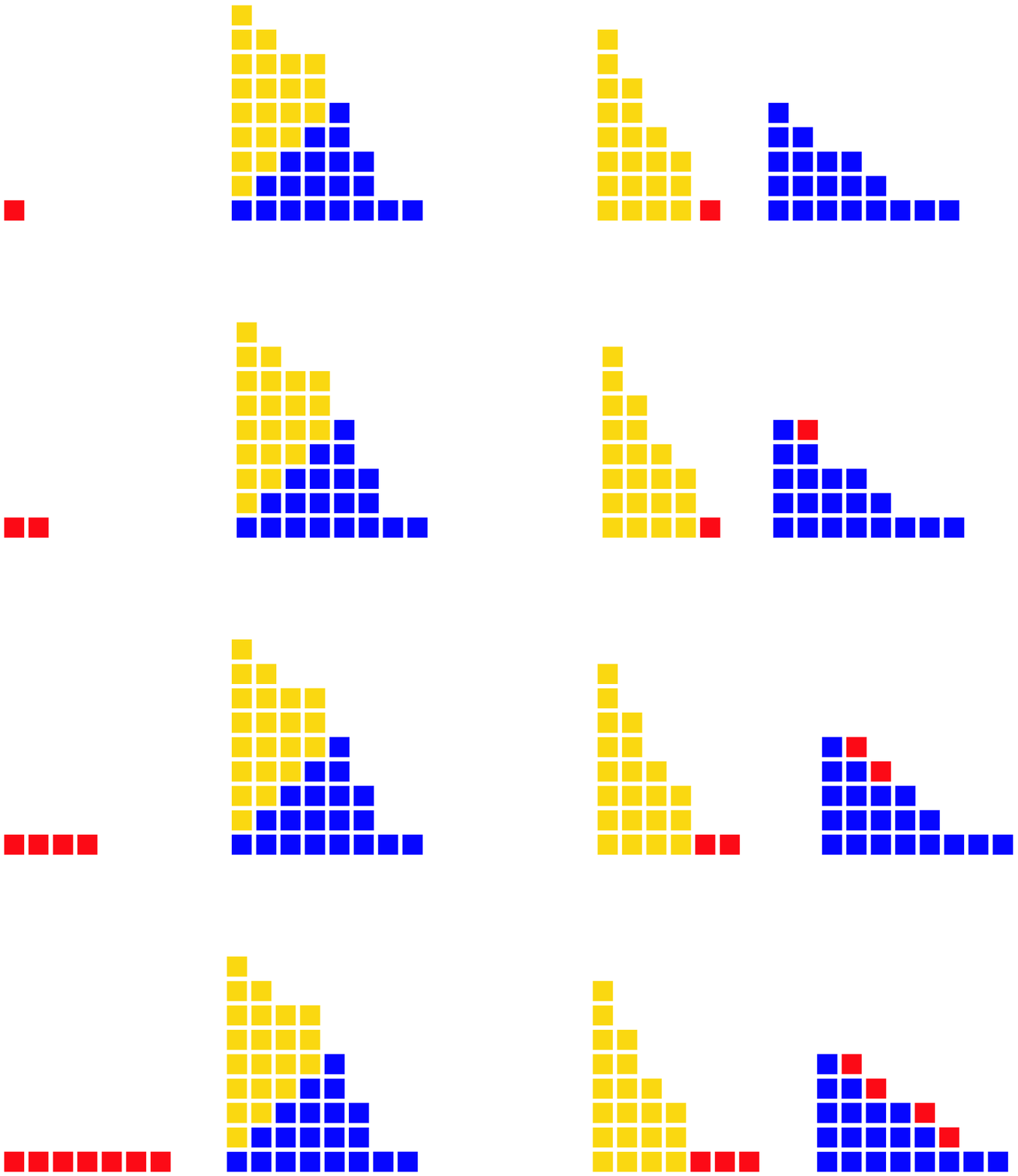}   
\caption{   
\label{recursif}  }  
\end{figure*}

 \noindent In Figure \ref{recursif} we illustrate the effect of the $*$-operation on pairs of the form $((n),\nu)$.

\noindent Given partitions $\mu$ and $\nu$, define the partition $\mu+\nu$ by 
   $$(\mu+\nu)_i:=\mu_i+\nu_i.$$
and set 
 $$\mu\cup \nu:=(\mu^{\prime} + \nu^{\prime})^{\prime}.$$
\noindent For example, 
if $\mu=33221$ and $\nu=531$, then 
$\mu \cup \nu = 53332211$ and $\mu+\nu=86321$.  
As usual, for $\mu$ and $\nu$ two partitions of $n$,  
$\mu$ is said to be {\em dominated} by $\nu$, in formula $\mu \preceq \nu$, 
 if for all $k \ge 1$: 
\[ 
\mu_1+ \mu_2+ \cdots+ \mu_k \leq \nu_1+ \nu_2+ \cdots+ \nu_k. 
\] 
Another remarkable property of the $*$-operation is that its image behaves nicely under the 
dominance order. More precisely,  
\begin{lemm}\label{dominance}   
 For any pair of partitions $(\mu,\nu)$, if $(\lambda,\rho)={(\mu,\nu)}^*$, then we have 
\begin{align} 
\mu \cup \nu & \succeq \lambda \cup  \rho,\qquad {\rm and\ equivalently} \label{union} \\ 
\mu + \nu & \preceq \lambda +  \rho.   \label{sum}   
\end{align} 
\end{lemm}   
\noindent
Observe that when $s_{\theta}$ appears in $s_{\mu}s_{\nu}$ with a nonzero coefficient, then  
\begin{equation*}\label{ineq} 
\mu\cup\nu \preceq \theta \preceq \mu + \nu, 
\end{equation*} 
thus (\ref{union}) and (\ref{sum}) imply that  
\[ \lambda\cup \rho\preceq \theta \preceq \la + \rho,\] 
which is compatible with Conjecture \ref{fflp}. 
 
\noindent Lemma \ref{dominance} immediately implies a statement very similar to that of Conjecture \ref{fflp}.  As is usual (See \cite{macdonald}), $h_\mu$ denotes the {\em complete homogeneous} symmetric function:
    $$h_\mu:=h_{\mu_1} h_{\mu_2} \cdots h_{\mu_k},$$
 with $h_a:=s_a$. 
 \begin{prop}\label{homogene}   
  For any pair of partitions $(\mu,\nu)$, if $(\lambda,\rho)={(\mu,\nu)}^*$, then
    $$h_{\lambda}h_{\rho}-h_\mu h_\nu$$ 
    is Schur-positive. 
 \end{prop} 
\noindent Recalling that $h_\mu h_\nu=h_{\mu\cup \nu}$, this follows from the fact that a difference of two homogeneous symmetric functions $h_{\alpha}-h_{\beta}$ is Schur-positive, if and only if $\alpha \preceq \beta$ (see \cite[Chapter 2]{sagan}).  A clear link between this proposition and Conjecture \ref{fflp} is established through the classical identity: 
   $$h_\alpha=s_\alpha + \sum_{\beta \succeq \alpha} K_{\beta\alpha}s_\beta,$$ 
where as usual $K_{\beta\alpha}$, the {\em Kostka} numbers,  count the number of semistandard tableaux of shape $\beta$ and type $\alpha$.  
 
The following result shows that the $*$-operation is also compatible with ``inclusion'' of partitions. Here, we say that $\alpha$ is {\em included} in $\mu$, if the diagram of $\alpha$ is included in the diagram of $\mu$.  We will simply write
    $$(\alpha,\beta)\subseteq (\mu,\nu),\qquad {\rm whenever}\qquad \alpha\subseteq \mu \quad {\rm and}\quad
         \beta\subseteq \nu.$$ 
\begin{lemm} \label{skew} Let $\alpha$, $\beta$,  $\mu$ and $\nu$ be partitions such that  $(\alpha,\beta) \subseteq (\mu,\nu)$. Then $\lambda(\alpha,\beta)\subseteq\lambda(\mu,\nu)$  and $\rho(\alpha,\beta)\subseteq\rho(\mu,\nu)$. 
\end{lemm} 
 
\noindent An immediate, but interesting, consequence of this lemma is the following observation. 
\begin{obs}\label{fixedshapes}
Let $(\alpha,\beta)$ and $(\gamma,\delta)$ be two fixed points of the $*$-operation such that  $(\alpha,\beta) \subseteq (\gamma, \delta)$. Writing simply $\lambda$ for $\lambda(\mu,\nu)$ and $\rho$ for $\rho(\mu,\nu)$, we see (using Lemma \ref{skew}) that
  $$(\alpha,\beta) \subseteq (\mu,\nu)  \subseteq (\gamma,\delta),$$
implies
 $$(\alpha,\beta) \subseteq (\lambda,\rho)  \subseteq (\gamma,\delta).$$
\end{obs}
\noindent As is underlined in \cite{fulton}, a pair of partitions $(\alpha,\beta)$ is a fixed point of the $*$-operation if and only if     
\begin{equation}\label{fixedpoints} 
\beta_1\geq \alpha_1 \geq \beta_2 \geq \alpha_2 \geq \cdots \geq \beta_n \geq \alpha_n.
\end{equation} 
Let us underline here that, for any $(\mu,\nu)$, it is easy to characterize the ``largest'' (resp. ``smallest'') fixed point contained in (resp. containing) the pair $(\mu,\nu)$. 
 We will see below how this observation can be used to link properties of $\lambda$ and $\rho$ to properties of $\mu$ and $\nu$.  Recall that an {\em horizontal strip} is a skew shape $\mu/\alpha$ with no two squares in the same column, and that a {\em ribbon} is a connected skew shape with no $2\times 2$ squares (see \cite[ Chapter 7]{stanley},  for more details). If we drop the condition of being connected in this last definition, we say that we have a {\em  weak ribbon}.  
 
Another striking consequence of Lemma \ref{skew} is that it allows a natural extension of the  $*$-operation to skew partitions. Denoting $(\mu,\nu)/(\alpha,\beta)$ the pair of skew shapes $(\mu/\alpha,\nu/\beta)$, we can simply define
\begin{equation}\label{extension}
   {(\mu/\alpha,\nu/\beta)}^*:={(\mu,\nu)}^*/{(\alpha,\beta)}^*.
\end{equation}
In other words, we have
\begin{equation}\label{lambda_extension} 
   \lambda(\mu/\alpha,\nu/\beta):=\lambda(\mu,\nu)/\lambda(\alpha,\beta), 
 \end{equation}
 and  
\begin{equation}\label{rho_extension} 
    \rho(\mu/\alpha,\nu/\beta):=\rho(\mu,\nu)/\rho(\alpha,\beta). 
   \end{equation}  
\noindent The  $*$-operation, or its extension as above,  preserves  (among others) the following families of pairs of (skew) shapes. 
\begin{prop} \label{other2} 
The  $*$-operation preserves the families of  
\begin{itemize} 
\item[\bf (1)] pairs of hooks; 
\item[\bf (2)] pairs of two-rows partitions;  
\item[\bf (3)] pairs of horizontal strips; 
\item[\bf (4)] pairs of weak ribbon. 
\end{itemize} 
\end{prop} 
 
\noindent Note that {\bf (1)} and {\bf (2)} follow directly from Observation \ref{fixedshapes}, and that the statements
{\bf (3)} and {\bf (4)} are made possible in view of our extension of the $*$-operation. 
 
\begin{figure*}[ht]  
\centering  
\begin{picture}(0,0)(0,0)
\put(111,24){$*$}
\put(105,20){$\longrightarrow$}
\end{picture} 
\includegraphics[width=80mm]{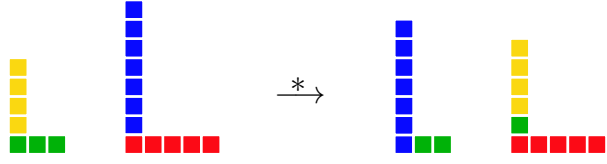}  
\caption{ The effect of the $*$-operation on hooks. 
\label{hooks}  } 
\end{figure*} 
 \noindent Results outlined in the sequel, and extensive computer experimentations suggests that we have the following extension of Conjecture \ref{fflp}. 
 
\begin{conj}\label{Bergeron} 
For any skew partitions $\mu/\alpha$ and $\nu/\beta$, if 
      $$(\lambda,\rho)={(\mu/\alpha,\nu/\beta)}^*,$$
then the symmetric function  
 \begin{equation}\label{conj_skew}
      s_{\lambda}s_{\rho}-s_{\mu/\alpha} s_{\nu/\beta}
\end{equation}
is Schur-positive. 
\end{conj} 
\noindent This has yet to be understood in geometrical terms. One should point out that there many skew shapes giving the same expression for the symmetric function $s_{\mu/\alpha} s_{\nu/\beta}$. The result of the $*$-operation is dependent on the particular choice of the skew-shape, so that there are many identities encoded in (\ref{conj_skew}). On the other hand, it is clear that Proposition \ref{cortranspose} extends to skew partitions. 
Our last combinatorial observation concerning the $*$-operation is the following. Let $\tau$ and $\nu$ be two fixed partitions, and consider all possible $\mu$'s such that $\rho(\mu,\nu)=\tau$. We claim that there is a minimal such $\mu$, if any, and we denote it $\theta(\tau,\nu)$. More precisely, we easily show that

\begin{prop}
 Given partitions $\tau$ and $\nu$, then for any $\mu$ such that $\rho(\mu,\nu)=\tau$, we must have
 \begin{equation}\label{rho_stable}
      \theta(\tau,\nu)\subseteq \mu.
  \end{equation}
 Furthermore, $\theta=\theta(\tau,\nu)$ is exactly the partition
      $$\theta=\tau_1^{b_1}\tau_2^{b_2-b_1}\tau_3^{b_3-b_2}\cdots,$$
with $b_j=\tau_j-\nu_j+j-1$.
\end{prop}
Condition (\ref{rho_stable}) can clearly be reformulated as
 \begin{equation}\label{rho_stable2}
    \mu_i\geq \tau_j,\qquad{\rm whenever}\quad  b_{j-1}<i\leq b_j,
 \end{equation}
 with $b_0:=0$.

\section{Main results}   
 
In this section we state our results concerning the validity of Conjecture \ref{fflp} for certain families of pairs, as well as its reduction to a finite number of tests for other families.  
We will show (in Section \ref{proof_general}) the following.
 
\begin{theo}\label{theo2rows}
Conjecture \ref{fflp} (or \ref{Bergeron}) holds
\begin{enumerate}
\item[\bf (1)] \label{theohooks} For  any pair $(\mu,\nu)$ of hook shapes. 
\item[\bf (2)] \label{twolines} For pairs of two-line and two-column partitions.
\item[\bf (3)] \label{skewhooks} For skew pairs the form $(\mu/\alpha,\nu/\beta)$, with all of $\mu$, $\nu$, $\alpha$ and $\beta$ are hooks.
\item[\bf (4)] \label{horiz} For skew pairs of the form $(0,\nu/\beta)$, with $\nu/\beta$ a weak ribbon.
\end{enumerate}
\end{theo} 

\noindent Other families for which we have partial results  correspond to Stembridge's (see \cite{stembridge}) classification of all  
multiplicity-free products of Schur functions. More precisely, these are all products of two Schur functions 
with Schur function expansion  having no coefficient larger then $1$. Thus, to show Conjecture \ref{fflp} in these cases, we need only show that the coefficient $c_{\lambda\rho}^{\theta}$ of $s_\theta$ in the product of $s_\lambda$ and $s_\rho$ is nonzero, whenever $c_{\mu\nu}^{\theta}=1$. 
 
 Stembridge uses the following notions for his presentation. A {\em rectangle} is a partition with at most one part size, i.e.,  empty, or of the form $(c^r)$ for suitable $c,r >0$; 
a  {\em fat hook} is  a partition with exactly two parts sizes, 
i.e., of the form $(b^r c^s)$ for suitable $b>c>0$; and 
a {\em near-rectangle}  is a fat hook such that it  
is possible to obtain a rectangle from it by  
deleting a single row or column. He shows that
the product $s_{\mu}s_{\nu}$ is multiplicity-free if and 
only if 
\begin{enumerate} 
\item[\bf (a)] $\mu$ or $\nu$ is a one-line rectangle, or 
\item[\bf (b)] $\mu$ is a two-line rectangle and $\nu$ is a 
fat hook or vice-versa, or 
\item[\bf (c)] $\mu$ is  rectangle and $\nu$ is a near rectangle  
 or vice-versa, or 
\item[\bf (d)]  $\mu$ and $\nu$ are both rectangles. 
\end{enumerate} 
 
Although we currently have proofs of the conjecture for cases {\bf (a)}  and {\bf (d)}  of Stembridge's pairs, proofs for cases {\bf (b)}  and {\bf (c)} are still in the process of being completed. Since all these share a common approach, we have decide to postpone their presentation to an upcoming paper.

 On another register, a careful study of the recursive construction of $\lambda(\mu,\nu)$ and $\rho(\mu,\nu)$ shows that,  
in a sense, Conjecture \ref{fflp} follows, under some conditions, from a finite number of cases when $\nu$ is fixed and $\mu$ becomes large.

More precisely, we obtain the result below. 
 As usual, the number of nonzero parts of $\mu$ is denoted by $\ell(\mu)$ and called  the {\em height} of $\mu$.

\begin{theo} \label{asympt} 
For any positive integer $p$, let $\nu$ be a fixed partition with at most $p$ parts, i.e. $\ell(\nu)\leq p$. Then, the validity of Conjecture \ref{fflp} for the infinite set of all pairs $(\mu,\nu)$, with $\ell(\mu) \leq p$, reduces to checking the validity of the conjecture for the finite set of pairs $(\alpha,\nu)$,  with $\alpha$ having at most $p$ parts, and largest part bounded as follows
\begin{equation}\label{borne}
   \alpha_1 \leq p\,(\nu_1+p).
\end{equation}
\end{theo}
\noindent
Theorem \ref{asympt} can also be generalized in a straightforward manner to the set of skew shapes pairs $(\mu/\alpha,\nu/\beta)$ of bounded height, with $\nu$ and $\alpha$ fixed.

\section{Proofs of the combinatorial properties}   

In the sequel, unless it is specifically mentioned, all partitions will be considered to have $n$ (possibly zero) parts.
We first observe that the set 
    $A_k(\mu,\nu):=\{j\ |\  \nu_j-j\geq \mu_k-k\},$
 appearing in (\ref{definition}), has to be of the form
    $A_k(\mu,\nu)=\{1,2,\ldots,a_k\}$
  for some $a_k=a_k(\mu,\nu)$, since $\nu_1-1>\nu_2-2>\ldots >\nu_{a_k}-a_k\geq\mu_k-k$, and thus 
  \begin{equation}\label{taille}
     a_k(\mu,\nu):=\#A_k(\mu,\nu).
  \end{equation}
In other words, 
    \begin{equation}\label{le_ak}
          \mu_k-k\leq \nu_m-m\qquad {\rm iff}\qquad 1\leq m\leq a_k.
      \end{equation}
Thus  definition (\ref{definition}) of $\lambda_k(\mu,\nu)$ can be reformulated as
 \begin{equation}\label{lambda}
     \lambda_k:=\mu_k-k+a_k(\mu,\nu).
\end{equation}
In the same spirit, we consider the set
   $B_j(\mu,\nu):=\{k\ |\   \mu_k-k> \nu_j-j\},$
which also has to be of the form $\{1,2,\ldots,b_j\}$, with
  \begin{equation}\label{tailleb}
     b_j(\mu,\nu)=\#B_j(\mu,\nu),
  \end{equation} 
  In other words, 
    \begin{equation}\label{le_bj}
          \nu_j-j< \mu_m-m\qquad {\rm iff}\qquad 1\leq m\leq b_j
      \end{equation}
and
 \begin{equation}\label{rho}
     \rho_j:=\nu_j-j+1+b_j(\mu,\nu).
\end{equation}

\begin{proof}[\bf Proof of Proposition \ref{recurrence}]\quad  To prove our recursive formula for the computation of the $*$-operation, we first analyze the case $\mu={\overrightarrow{\alpha}}^i$. As we have already mentioned, this means that
$\mu_k=\alpha_k$ for all $k\neq i$, and $\mu_i=\alpha_i +1$. 
Let $(\lambda,\rho)={(\alpha,\nu)}^*$. Now, suppose that there exists a $ j \in \{1,\ldots,n\}$ such that
\[\nu_j - j = \alpha_i -i .\]
This implies that $A_k(\mu,\nu)=A_k(\alpha,\nu)$ for all $k \neq i$
and that
\[A_i(\alpha,\nu)=\{1,2,\ldots,j\},\;\; {\rm and} \;\; A_i(\mu,\nu)=\{1,2,\ldots,j-1\}.\]
It follows that $\lambda_k(\mu,\nu)=\lambda_k$ for all $k\neq i$, and that 
\begin{eqnarray*} \lambda_i(\mu,\nu) & = & \mu_i -i + a_i(\mu,\nu)\\
                 & = & \alpha_i +1 -i + a_i(\alpha,\nu)-1
                  =  \lambda_i
\end{eqnarray*}
Hence $\lambda(\mu,\nu)=\lambda$.

\noindent On the other hand, we clearly have $B_k(\alpha,\nu)=B_k(\mu,\nu)$ for all $k \neq j$, and
\[B_j(\alpha,\nu)=\{1,2,\ldots,i-1\} \;\; {\rm and} \;\; B_j(\mu,\nu)=\{1,2,\ldots,i\}.\]
Hence $\rho_k(\mu,\nu)=\rho_k$ for all $k\neq i$ and,
\begin{eqnarray*} \rho_j(\mu,\nu) & = & \nu_j -j +1+ b_j(\mu,\nu)\\
                 & = & \alpha_i -i +1 + i
                  =  \mu_i.
\end{eqnarray*}
Since $\rho_j=\mu_i-1$ we conclude that $\rho(\mu,\nu)={\overrightarrow{\rho}}^j$, and this settles the first case of (\ref{gauche}).

If no $ j \in \{1,\ldots,n\} $ is such that $\nu_j - j = \alpha_i -i $, we have the equalities 
  $$A_k(\mu,\nu)=A_k(\alpha,\nu)\qquad {\rm and} \qquad B_k(\mu,\nu)=B_k(\alpha,\nu)$$ 
  for all $k \in \{1,\ldots,n\}$. It follows $\lambda_k(\mu,\nu)=\lambda_k$ for all $k \neq i$ and
\begin{eqnarray*}
\lambda_i(\mu,\nu) = \alpha_i+1 -i + a_i(\alpha,\nu) =  \lambda_i +1,
\end{eqnarray*}
so $\lambda(\mu,\nu)={\overrightarrow{\lambda}}^i $. It easily follows that in this case  $\rho(\mu,\nu)=\rho$ and this
concludes the second case.
The part (\ref{droite}) of the proposition is shown in a similar manner.
\end{proof}


\noindent As announced in (\ref{transpose}), we have the following lemma.
\begin{lemm}\label{conjugate}
For any pair of partitions $(\mu, \nu)$, we have
  $$  {(\mu,\nu)}^*=(\lambda,\rho)\qquad{\rm iff}\qquad {(\nu',\mu')}^*=(\lambda',\rho').$$
\end{lemm}
\begin{proof}[{\bf Proof.}]\quad  We proceed by induction on $|\mu|+|\nu|$. The lemma obviously holds when
$\mu=\nu=0$, so suppose that (\ref{transpose})  holds for all $(\alpha, \nu)$ with 
$\alpha \subseteq \mu$. That is,
    $$\lambda(\alpha, \nu) = \lambda'(\nu', \alpha') \qquad     
               {\rm and } \qquad
         \rho(\alpha, \nu) = \rho'(\nu', \alpha'). $$
In view of the recursive description of the $*$-operation, 
it is easy to verify that we need only show that the ``if'' part of   (\ref{gauche}) applies to the pair $(\mu, \nu)$ if and only
the ``otherwise'' part of (\ref{droite}) applies to the pair $(\nu', \mu') $. Let us suppose that $\mu=\overrightarrow{\alpha}^i$. We want to show that
 $$\nu_j -j = \alpha_i - i$$
if and only if there is no $k$ such that $\nu'_k - k = \alpha'_{\mu_i} - \mu_i +1$. 
Assume
that there exists a $j$ such that $\alpha_i-i = \nu_j -j$. (See Figure \ref{corner}.)
\begin{figure}[ht]  
\centering 
\vskip10pt
\begin{picture}(10,20)(00,0)
\put(50,250){$\mu={\overrightarrow{\alpha}}^i$}
\put(180,260){$\mu'=\alpha'\col_{i}$}
\put(0,205){$i{\scriptstyle\rightarrow}$}
\put(148,167){$\mu_i{\scriptstyle\rightarrow}$}
\put(30,133){$\uparrow\atop\displaystyle\mu_i$}
\put(225,133){$\uparrow\atop \displaystyle i$}
\put(60,90){$\nu$}
\put(210,90){$\nu'$}
\put(-2,66){$j {\scriptstyle\rightarrow}$}
\put(148,18){$\nu_j{\scriptstyle\rightarrow}$}
\put(28,-15){$\uparrow\atop\displaystyle\nu_j$}
\put(235,-15){$\uparrow\atop \displaystyle j$}
\end{picture}
\includegraphics[width=90mm,height=88mm]{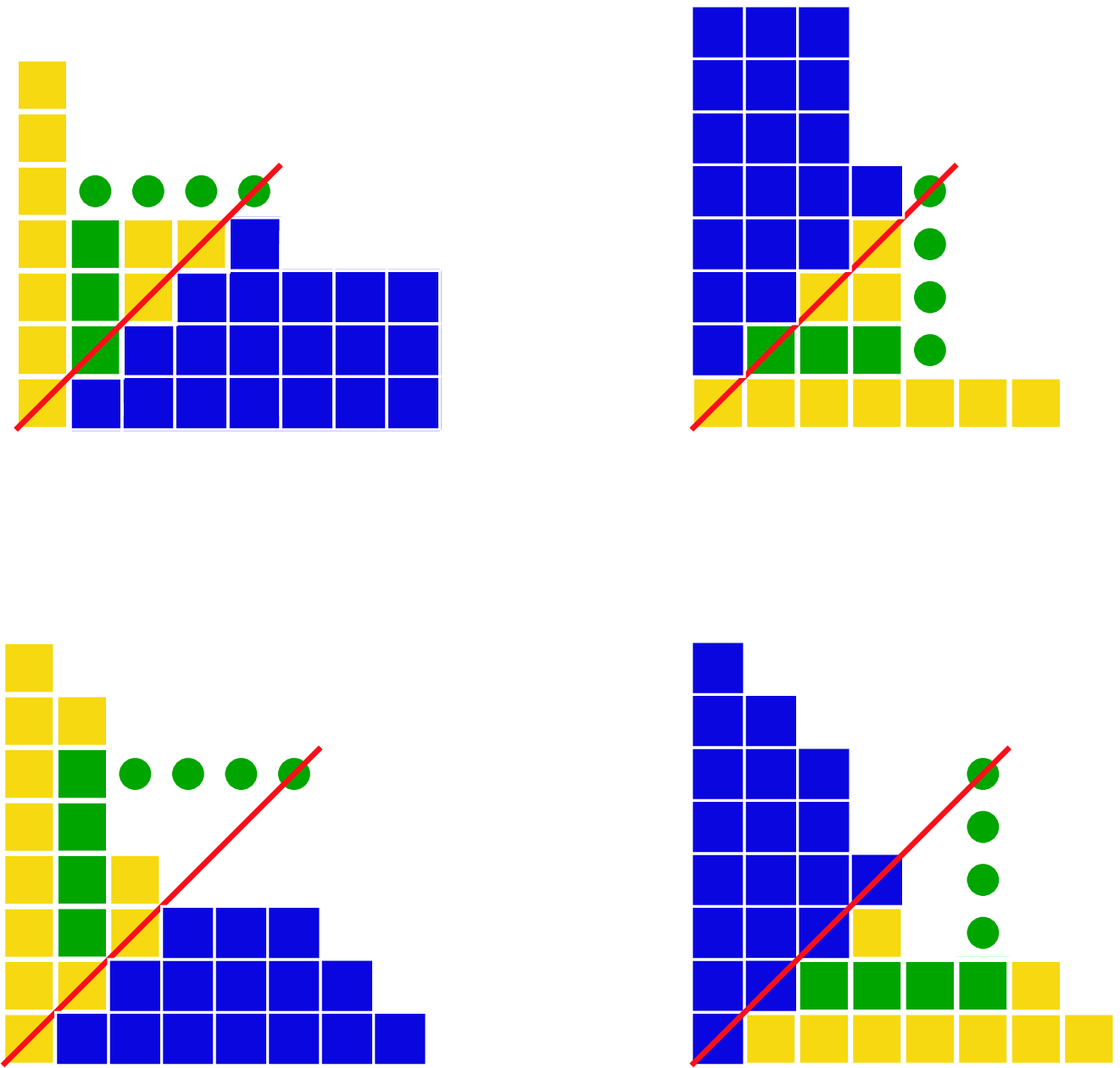}
\bigskip\caption{}\label{corner}
\end{figure}
Then, by (\ref{gauche}) we have 
 $$\lambda(\alpha, \nu) = \lambda(\mu,\nu)=\lambda'(\nu', \alpha')= \lambda'(\nu', \mu'),$$
and
   $$ \rho(\mu,\nu)={\overrightarrow{\rho(\alpha, \nu)}}^j.$$
 We thus need to check that
   $$\rho'(\nu', \mu')=\rho'(\nu', \alpha') \col_{j} .$$
Considering the pair $(\nu', \mu')$, we have $\mu'={\overrightarrow{\alpha' }}^{\mu_i} $.
Thus, the new cell is added  to $\alpha'$ in row $\mu_i$, which is of
length $i-1$, i.e.: $\alpha'_{\mu_i}=i-1$. We claim that there does not exist a
 row $k$ of $\nu'$ such that 
 $$\nu_k'-k = (i-1) - (\alpha_i +1) + 1 = i -\alpha_i-1.$$ 
In fact, if $k=\nu_j$ the difference   $\nu_k'-k$ is strictly 
bigger than $i -\alpha_i-1$, since 
    $$\nu'_k - k \ge j - \nu_j = i - \alpha_i.$$
On the other hand, if $k=\nu_j+1$, the difference   $\nu_k'-k$ has to be strictly 
smaller than $i -\alpha_i-1$, since $(\nu_{j+1}, j) \notin \nu$, and hence
$(j, \nu_{j+1}) \notin \nu'$.  Then, 
  $$\nu'_{\nu_j + 1} - ( \nu_j + 1) < j - (\nu_j + 1) = i - \alpha_i -1.$$
\end{proof}
\noindent In preparation for the proof of Lemma \ref{dominance}, let us prove the following.
\begin{lemm}\label{dominance_sum}The following two statements are equivalent
\begin{itemize}
\item[\bf (1)] For all $(\mu,\nu)$ pair of partitions $\mu \cup \nu \succeq \la \cup \rho$.
\item[\bf (2)] For all $(\mu,\nu)$ pair of partitions $\mu + \nu \preceq \la + \rho$.\end{itemize}
\end{lemm}
\begin{proof}[{\bf Proof.}]\quad 
Assuming {\bf(1)}, we have
\begin{eqnarray*}
   \mu+\nu&=& (\mu'\cup \nu')'\\
                 & \preceq & (\lambda(\nu',\mu')\cup\rho(\nu',\mu'))'\\
                 &=& \lambda'(\nu',\mu')+\rho'(\nu',\mu')\\
                 &=& \lambda(\mu,\nu)+\rho(\mu,\nu).
\end{eqnarray*}
A similar computation shows the reverse statement.
\end{proof}


\begin{proof}[\bf Proof of Lemma \ref{dominance}] (P. McNamara)\quad 
Let $(\lambda,\rho)=(\mu,\nu)^*$, for $\mu$ and $\nu$ partitions with $n$ (possibly zero) parts. From the above observation, it is sufficient to show that 
  $$   \lambda+\rho\succeq    \mu+\nu ,$$
and hence that, for all $i$,  
  \begin{equation}\label{a_montrer}
   \sum_{j=1}^i (\lambda_j+\rho_j)\geq \sum_{j=1}^i (\mu_j+\nu_j).
   \end{equation}
Definitions (\ref{lambda}) and (\ref{rho}) give $\lambda_j=\mu_j-j+a_j$ and $\rho_j=\nu_j-j+1+b_j$,  so that (\ref{a_montrer}) becomes
  \begin{equation}\label{a_montrer2}
   \sum_{j=1}^i (\mu_j+\nu_j+a_j+b_j-(2\,j-1))\geq \sum_{j=1}^i (\mu_j+\nu_j),
   \end{equation}  
which is equivalent to 
\begin{equation}\label{quadrato}
\sum_{j=1}^i (a_j+ b_j)\ge i^2.
\end{equation}
But the definitions of $a_j$ and $b_j$, can clearly be reformulated as
   $$a_j=\#\{ (k,j)\ |\ \nu_k-k\geq \mu_j-j\ \}\qquad{\rm and}\qquad
       b_k=\#\{ (k,j)\ |\ \nu_k-k< \mu_j-j\ \}.$$
Thus the left hand side of (\ref{quadrato}) is the cardinality of the set
   $$\{ (k,j)\ |\ 1 \leq k\leq i, \ 1\leq j\leq n \quad {\rm or}\quad 1\leq k\leq n,\ 1\leq j\leq i\ \}$$
which clearly contains the set $\{1,\ldots,i\}\times \{1,\ldots,i\}$, hence the inequality.
   \end{proof}
\begin{proof}[\bf Proof of Lemma \ref{skew}]\quad  We now want to show that $(\alpha,\beta)\subseteq (\mu,\nu)$ implies that $(\alpha,\beta)^*\subseteq (\mu,\nu)^*$, for partitions with $n$ parts. For this, we need to prove that, for every $1\le k \le n$, 
      $$\lambda_k(\alpha,\beta)\le \lambda_k(\mu,\nu)\qquad {\rm and}\qquad
          \rho_k(\alpha,\beta)\le \rho_k(\mu,\nu).
       $$
Recall that by definition 
   $$ \lambda_k(\alpha,\beta)=\alpha_k -k + a_k(\alpha,\beta) \;\;{\rm and} \;\;  
   \lambda_k(\mu,\nu)=\mu_k -k + a_k(\mu,\nu).$$
Thus  there is nothing to show when $a_k(\alpha,\beta)\leq a_k(\mu,\nu)$, since by hypothesis $\alpha_k\leq \mu_k$. Otherwise, in the case when 
    $$a_k(\alpha,\beta)- a_k(\mu,\nu)>0,$$
the inequality we wish to show is clearly equivalent to
\begin{equation}\label{toprove}
   (\mu_k-k) -(\alpha_k-k) \geq a_k(\alpha,\beta) - a_k(\mu,\nu).
\end{equation}
Using (\ref{le_ak}) and the inclusion $\beta\subseteq\nu$,  we see that for all $m$, with $a_k(\mu,\nu)<m\leq a_k(\alpha,\beta)$, we must  have      
\begin{equation}
    \label{inegalite}
      \alpha_k-k\leq \beta_m-m\leq \nu_m-m<\mu_k-k.
\end{equation}
 Observing that all the $(\beta_m-m)$'s appearing in (\ref{inegalite}) are pairwise distinct, we have made evident inequality (\ref{toprove}) since  at least $a_k(\alpha,\beta) - a_k(\mu,\nu)$ integers  separate $\alpha_k-k$ from $\mu_k-k$.
 \end{proof}


\begin{proof}[{\bf Proof of Proposition \ref{other2}}]\quad  
Among the families of pairs stated to be preserved by the $*$-operation, we have already shown cases {\bf (1)} and {\bf (2)}. The proofs of the other two claims are as follows.
\item[{\bf (3)}] Recall that $\mu/\alpha$ and $\nu / \beta$ are horizontal strips if and only if, for all $1\leq k<n$
 $$\alpha_k \geq \mu_{k+1}\qquad{\rm and}\qquad \beta_k \geq \nu_{k+1}.$$
To show that $(\mu/\alpha,\nu/\beta)^*$ is also an horizontal strip, we need to prove that $\lambda_k(\alpha,\beta)\geq \lambda_{k+1}(\mu,\nu)$ and $\rho_k(\alpha,\beta)\geq \rho_{k+1}(\mu,\nu)$. The approach is similar to that of Lemma \ref{skew}. Once again, by definition, we have
$$\rho_k(\alpha,\beta)=\beta_k -k+1 +b_k(\alpha,\beta)
       \quad{\rm and} \quad\rho_{k+1}(\mu,\nu)=\nu_{k+1} -k +b_{k+1}(\mu,\nu).$$
If $b_k(\alpha,\beta)\geq b_{k+1}(\mu,\nu)$ there is nothing to prove. Otherwise, 
the inequality we want to prove is clearly equivalent to 
\begin{equation}\label{claime}
   b_{k+1}(\mu,\nu) - b_k(\alpha,\beta)\leq \beta_k-\nu_{k+1}+1.
\end{equation} 
Using (\ref{le_bj}), when $b_k(\alpha,\beta)<m\leq b_{k+1}(\mu,\nu)$, we must have $\alpha_m-m\leq \beta_k-k$, and   $\nu_{k+1}-(k+1)<\mu_m-m$. Hence
\begin{equation*}
     \nu_{k+1}-(k+1)<\mu_{m+1}-(m+1)<\alpha_m-m\leq \beta_k -k,
\end{equation*}
since by hypothesis $\mu_{m+1}\leq \alpha_m$.
This shows that there are at least $b_{k+1}(\mu,\nu)-b_k(\alpha,\beta)-1$ distinct integers separating $\nu_{k+1}-(k+1)$ from $\beta_k -k$, thus (\ref{claime}) follows. Similarly we get the other inequality.

\item[{\bf (4)}] The statement that  $\mu/\al$ and $\nu/\beta$ are two weak ribbons is equivalent to saying that
\begin{equation}\label{rib}   
\mu_{k+1}\leq \alpha_k+1 \qquad {\rm and}\qquad \nu_{k+1} \leq \beta_k+1,
\end{equation}
for all $1\leq k<n$. 
If $\mu_{k+1}\leq \alpha_k$ then the result follows by part {\bf (3)}. So suppose $\mu_{k+1}=\alpha_k+1$. We have 
\begin{equation}\label{stesso}
\mu_{k+1}-(k+1)=\alpha_k-k.
\end{equation}
Thus, showing $\la_{k+1}(\mu,\nu)\leq \la_k(\al,\be)+1$ is equivalent to showing 
\[a_{k+1}(\mu,\nu)\leq a_k(\al,\be)+1.\]
But this easily follows by the definitions of $A_k(\al,\be)$, $A_{k+1}(\mu,\nu)$, and from (\ref{rib}) and (\ref{stesso}). The proof of the other inequality is similar. This last case concludes the proof.
\end{proof}

\section{Extension of the $*$-operation to tableaux}\label{ext_tableau}
Since we are using the french notation for partitions, standard tableaux have increasing entries  along rows from left to right, and increasing entries along columns from bottom to top. As usual, a semistandard tableau is one in which we relax the requirement along rows to weakly increasing. The {\em reading word} of a tableau is obtained by reading the entries of the tableau starting with the top row, from left to right, and going down  the rows.  For instance, the reading word of
$$  \begin{matrix}\Young{\Box{5}&\Box{5}\cr\Box{2}&\Box{3}&\Box{4}\cr\Box{1}&\Box{2}&\Box{3}&\Box{3}\cr} \end{matrix}$$
is $552341233$. It is  well known that semistandard tableaux correspond to a chain in the Young lattice, $ 0 \subseteq \mu_1 \subseteq \mu_2 \subseteq \ldots\subseteq \mu_k=\mu$, such that $\mu_{i+1}/\mu_i$ is a horizontal strip.
   
For instance the chain associated to the tableau above is
    $$0\subseteq 1\subseteq 21\subseteq 42\subseteq 43\subseteq 432,$$
and clearly standard tableaux correspond to maximal chains. The shape of a tableau is the final partition in the corresponding chain. All of these notions extend to skew shapes. In particular, a semistandard tableau of skew shape $\lambda/\mu$ is a chain starting at shape $\mu$ and ending at shape $\lambda$. For this to be possible, we clearly need $\mu\subseteq \lambda$. We sometimes say that a semistandard tableau, of shape $\lambda/\mu$, is a {\em filling} of $\lambda/\mu$. The {\em natural} filling of a partition $\mu=(\mu_1,\mu_2,\ldots,\mu_k)$, is the semistandard tableau corresponding to the chain
   $$0\subseteq(\mu_1)\subseteq(\mu_1,\mu_2)\subseteq(\mu_1,\mu_2,\mu_3)
          \subseteq\ldots \subseteq(\mu_1,\mu_2,\ldots,\mu_k)$$
Thus, each cell is {\em filled} by the number of the row it lies in.
The {\em type} of a semistandard (possibly skew shaped) tableau $t$ is the sequence $(m_1,m_2,\ldots)$ of multiplicities of its entries. This is to say that $m_i=m_i(t)$ is the number of entries that are equal to $i$. When $m_1\geq m_2\geq \ldots$, this type can be identified with a partition. The natural filling of $\mu$ is the only semistandard tableau of shape $\mu$ that also has type $\mu$.

\begin{figure*}[ht] 
\centering  
\includegraphics[width=120mm]{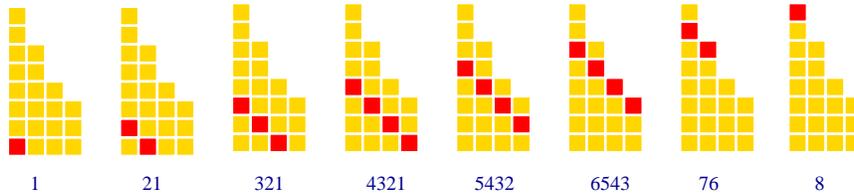}  \vskip-20pt 
\caption{ $\Delta(\mu)$: the  diagonal reading of the natural filling of $\mu$. 
\label{grow}  } 
\end{figure*} 
\noindent The {\em diagonal reading}, $\Delta(t)$, of  a tableau $t$, of shape $\mu$,  is  obtained by recording the entries of $t$ following the diagonals $x+y=k$ in the partition $\mu$, from left to right, and from top to bottom, for $k=1, 2, \ldots ,\ell(\mu)$. We simply denote $\Delta(\mu)$ the diagonal reading of the natural filling of $\mu$.
For    $\mu=44432211$, we have    
    $$\Delta(\mu)=1\, 21\, 321\, 4321\, 5432\, 6543\, 76\, 8$$
   as is illustrated in Figure \ref{grow}.  

To describe more consequences of the properties of $*$, we consider the {\em double Young lattice}, $\D$, which is just the direct product of two copies of the usual Young lattice. The double Young lattice already plays an explicit role in \cite{fulton}, see also \cite{fomsta}.

There is a natural  grading for $\D$ given by $(\mu,\nu)\mapsto |\mu|+|\nu|$. A {\em standard (tableau) pair} of {\em shape} $(\mu,\nu)$ is a maximal chain in this graded poset that starts at $(0,0)$  and ends at $(\mu,\nu)$. For example, we have 
\begin{equation}\label{standpair}
 (0,0)\subseteq (0,1)\subseteq (0,2)\subseteq(1,2)\subseteq (1\,1,2)\subseteq (2\,1,2)\subseteq (2\,1,3)
 \end{equation}
As in the usual case, such a chain can be identified with a pair $(t,r)$ of standard tableaux, of respective shapes $\mu$ and $\nu$, with non-repeated entries from the the set $\{1,2,\ldots,n\}$, $n=|\mu|+|\nu|$. The number $f_{(\mu,\nu)}$ of standard pairs of shape $(\mu,\nu)$ is thus  
 \begin{equation}\label{nombre_tableaux} 
      f_{(\mu,\nu)}={|\mu|+|\nu|\choose |\mu|}f_\mu\,f_\nu 
  \end{equation} 
  where $f_\mu$ and $f_\nu$ are both given by the usual hook formula. 
   In terms of tableaux, the standard pair (\ref{standpair}) corresponds to: 
     $$\begin{pmatrix} 
               \Young{\Box{4}\cr\Box{3}&\Box{5}\cr} &,& \Young{\Box{1}&\Box{2}&\Box{6}\cr} 
               \end{pmatrix}.$$ 
The double Young lattice occurs naturally in the study of representations of the hyperoctahedral groups. This suggests that there might be a link between that subject and the study of properties of the transformation $*$.    

A {\em semistandard pair} is a chain 
     $$(0,0)=\pi_0\subseteq\pi_1\subseteq \cdots \subseteq \pi_k=(\mu,\nu)$$ 
  in   $\D$, such that $\pi_{j+1}/\pi_j$ is an horizontal strip pair for each $1\leq j\leq k-1$. 
  For example, the pair of semistandard tableaux 
    $$\begin{pmatrix} 
               \Young{\Box{3}\cr\Box{2}&\Box{3}&\Box{3}\cr} &,& \Young{\Box{2}&\Box{3}\cr\Box{1}&\Box{1}&\Box{3}\cr} 
               \end{pmatrix},$$ 
 corresponds to the path
   $$(0,0)\subseteq (0,2)\subseteq (1,21)\subseteq (31,32).$$
 It follows from Proposition \ref{other2} that 
\begin{lemm} \label{tableaux}   
  The function $*:\D\longrightarrow\D$ is  an increasing transformation that preserves both standard and semistandard pairs. 
  \end{lemm}  
We can thus extend the $*$-operation to semistandard (and standard) pairs.
For example, 
    $$\begin{pmatrix} 
\Young{\Box{26}\cr\Box{22}&\Box{23}&\Box{24}\cr\Box{16}&\Box{17}&\Box{18}&\Box{19}&\Box{20}\cr\Box{9}&\Box{10}&\Box{11}&\Box{12}&\Box{13}\cr}&\!\!,\!\!&  
\Young{\Box{25}\cr\Box{21}\cr\Box{14}&\Box{15}\cr\Box{1}&\Box{2}&\Box{3}&\Box{4}&\Box{5}&\Box{6}&\Box{7}&\Box{8}\cr}               
 \end{pmatrix}^{\Large\! *}=\begin{pmatrix} 
\Young{\Box{26}\cr\Box{22}&\Box{24}\cr\Box{16}&\Box{17}&\Box{19}&\Box{20}\cr\Box{9}&\Box{10}&\Box{11}&\Box{12}&\Box{13}\cr} 
&\!\!,\!\!&  
\Young{\Box{25}\cr\Box{21}&\Box{23}\cr\Box{14}&\Box{15}&\Box{18}\cr\Box{1}&\Box{2}&\Box{3}&\Box{4}&\Box{5}&\Box{6}&\Box{7}&\Box{8}\cr} 
 \end{pmatrix}$$ 
We emphasize that the resulting filling of $(\lambda,\rho)$ heavily depends on the particular filling of $(\mu,\nu)$ that has been chosen as the input standard pair. Fixed points, for standard pairs, are easily characterized as follows.
 
 \begin{lemm} \label{fixed_tableaux}   
  A standard pair $(t,r)$, of shape $(\mu,\nu)$, is fixed point of the $*$-operation, if and only if $(\mu,\nu)$ is fixed, and the tableau, obtained by alternating rows of $r$ and rows of $t$,  is standard.
\end{lemm}
\noindent
Recall that, if the pair  $(\mu, \nu)$ is a fixed point, then (\ref{fixedpoints}) implies that the alternating lengths of the rows are in decreasing  order.

\section{Background on Littlewood-Richardson coefficients}
In order to prove that Conjecture \ref{fflp}, or our extension of it, holds for some given pairs, we clearly need one of the many classical descriptions of Littlewood-Richardson coefficients. From a broader perspective, let us briefly recall some classical facts about these coefficients (see \cite{macdonald} or \cite{stanley} for more details). For $\mu$ and $\nu$ two partitions, and $\theta$ such that $|\theta|=|\mu|+|\nu|$, the coefficient $c_{\mu\nu}^\theta$ of $s_\theta$ in $s_\mu s_\nu$ is clearly given as
 \begin{eqnarray}\label{scalar}
   c_{\mu\nu}^\theta&=&\langle s_\mu s_\nu,s_\theta\rangle\\
         &=&\langle s_{\theta/\mu},s_\nu\rangle,
\end{eqnarray}
where $\langle-,-\rangle$ denotes the usual scalar product on symmetric function, for which Schur functions are orthonormal. 

The following is the explicit formulation of the Littlewood-Richardson rule that we are going to use to compute the $c_{\mu\,\nu}^\theta$'s. In order to state it, let us recall some terminology.  A {\em lattice permutation} is a sequence of positive integers $a_1a_2\cdots a_n$ such that in any initial factor $a_1a_2\cdots a_j$ the number of $i$'s is at least as great as the number of $i+1$'s, for all $i$. The type of a lattice permutation is (naturally) the sequence of multiplicities of the integers $1,2,\ldots$ that appear in it. Note that $\Delta(\mu)$ (the diagonal reading of the natural filling of $\mu$) is always a lattice permutation  of type $\mu$.  The {\em reverse reading word} of a tableau, is the reading word of a tableau, read backwards. 
A proof of the following assertion can be found in \cite{stanley}.

\begin{enumerate}
\item[] {\em
{\bf Littlewood-Richardson Rule.}\quad The Littlewood-Richardson coefficient $c_{\mu\,\nu}^\theta$ is equal to the number of semistandard tableaux of shape $\theta /\nu$ and type $\mu$ whose reverse reading word is a lattice permutation. }
\end{enumerate}

\noindent When a semistandard tableaux of shape $\theta/\nu$ has a lattice permutation as its reverse reading word, we say that it is a  LR-{\em filling} of shape $\theta/\nu$. 

\begin{center}
\begin{figure*}[ht]
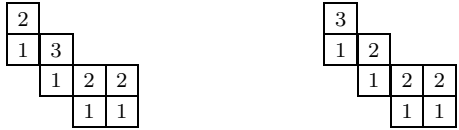
 
  $$\begin{matrix}\Young{\Box{2}\cr
                    \Box{1}&\Box{3}\cr
                    \Blk&\Box{1}&\Box{2}&\Box{2}\cr
                    \Blk&\Blk&\Box{1}&\Box{1}\cr}
           \qquad\qquad\qquad
      \Young{\Box{3}\cr
                    \Box{1}&\Box{2}\cr
                    \Blk&\Box{1}&\Box{2}&\Box{2}\cr
                    \Blk&\Blk&\Box{1}&\Box{1}\cr}\end{matrix}$$
\caption{The two LR-fillings of $4421/21$ of type $431$.}
\label{2tableaux}  
\end{figure*}
\end{center}

\noindent For $\theta=4421$, $\nu=21$ and $\mu=431$, we have $c_{\mu\,\nu}^\theta=2$ since there are exactly two LR-fillings of $\theta/\nu$.  These are described in Figure \ref{2tableaux}. The two corresponding reverse reading words $11221312$ and $11221213$. They are clearly lattice permutations of type $\mu$.

\section{Proof of special instances}\label{proof_general}

In this section we show that Conjecture \ref{fflp} holds for   pairs of hooks, pairs of two-row (or two-columns) shapes, and in a special case corresponding to our generalization of the conjecture to skew partitions. We first prove that Conjecture \ref{fflp} holds when one of the partitions is empty.

\begin{figure*}[ht] 
 \begin{picture}(0,0)(0,0)
 \put(20,-10){$\overline{\nu}$}
\put(60,30){\Huge$ \nearrow$}
\put(150,-10){$\nu$}
\put(208,73){$\circlearrowleft$}
\put(200,60){\Huge$ \rightarrow$}
\put(280,-10){$\nu/\underline{\nu}$}
\end{picture}\bigskip\bigskip
\centering  
\includegraphics[width=120mm]{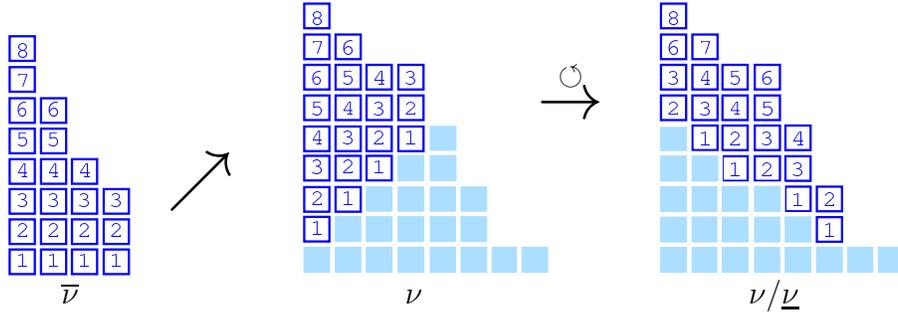} \vskip-15pt 
\caption{The LR-filling of $\nu / \underline{\nu}$ with reverse reading word  $\Delta(\overline{\nu})$. 
\label{filling}  } 
\end{figure*}

\begin{lemm}
  For any partition $\nu$, setting $\overline{\nu}:=\lambda(0,\nu)$ and $\underline{\nu}:=\rho(0,\nu)$, $\Delta(\overline{\nu})$ is the reverse reading word of a LR-filling of $\nu/\underline{\nu}$ of type $\overline{\nu}$.
\end{lemm}
\begin{proof}[{\bf Proof.}]\quad 
We show that   $\Delta(\overline{\nu})$ encodes a LR-filling of $\nu  / \underline{\nu}$ of type $\overline{\nu}$.  To this end, we proceed as follows. We ``slide'' the natural filling of $\overline{\nu}$  up the columns of $\nu$. This gives a partial filling of $\nu$ with empty cells for the portion of $\nu$ that corresponds to $\underline{\nu}$. We will suppose that these empty cells are filled with zeros.  We then sort each row in increasing order to get a filling of the skew shape $\nu/\overline{\nu}$.  By construction, we obtain a filling of $\nu  / \underline{\nu}$ whose reverse reading word is the lattice permutation $\Delta(\overline{\nu})$. An example is given in Figure \ref{filling}.

To show the lemma, we need only show that the resulting tableau is  semistandard. We already have strict increase along rows, so we need only check that this is also true along columns. By construction, the right-most entry in the $(k+1)^{\rm th}$-row of the final filling of $\nu/\underline{\nu}$ is $k$. Since the integers in a row are consecutive by construction, the difference between two entries in the same column of $\nu  / \underline{\nu}$, one in  the $i^{\rm th}$ row and the other in  the $(i+1)^{\rm th}$  row, has to be equal to $\nu_i - \nu_{i+1} +1$, which is larger than zero.
\end{proof}
\noindent It immediately follows that
\begin{coro}\label{empty}
For any partition $\nu$ the difference
    $$s_{\overline{\nu}} s_{\underline{\nu}}-s_\nu$$
 is Schur positive. Thus, recalling that $s_0=1$, Conjecture \ref{fflp} holds  for pairs of the form $(0,\nu)$.   
\end{coro} 

 \begin{figure} [ht] 
  \begin{center}
 \begin{picture}(0,0)(0,0)
\put(6,212){$\underbrace{\hskip22pt}_{\displaystyle a}$}
\put(176,212){$\underbrace{\hskip22pt}_{\displaystyle a}$}
\put(-13,234){$b\left\{\rule{0cm}{19pt}\right.$}
\put(158,240){$d\left\{\rule{0cm}{33pt}\right.$}
\put(56,212){$\underbrace{\hskip42pt}_{\displaystyle c}$}
\put(228,212){$\underbrace{\hskip42pt}_{\displaystyle c}$}
\put(37,248){$d\left\{\rule{0cm}{33pt}\right.$}
\put(210,243){$b\left\{\rule{0cm}{19pt}\right.$}
\put(126,239){$*$}
\put(120,235){$\longrightarrow$}
\put(70,183){{\bf (a)} Case $a\leq c$ and $b\leq d$}
\put(6,110){$\underbrace{\hskip22pt}_{\displaystyle a}$}
\put(176,110){$\underbrace{\hskip22pt}_{\displaystyle a}$}
\put(-13,146){$b\left\{\rule{0cm}{34pt}\right.$}
\put(158,132){$d\left\{\rule{0cm}{19pt}\right.$}
\put(56,110){$\underbrace{\hskip42pt}_{\displaystyle c}$}
\put(228,110){$\underbrace{\hskip42pt}_{\displaystyle c}$}
\put(37,132){$d\left\{\rule{0cm}{19pt}\right.$}
\put(210,146){$b\left\{\rule{0cm}{34pt}\right.$}
\put(126,137){$*$}
\put(120,133){$\longrightarrow$}
\put(70,81){{\bf (b)} Case $a\leq c$ and $b> d$.}
\put(5,0){$\underbrace{\hskip42pt}_{\displaystyle a}$}
\put(176,0){$\underbrace{\hskip33pt}_{\displaystyle a-1}$}
\put(-13,24){$b\left\{\rule{0cm}{19pt}\right.$}
\put(158,28){$d\left\{\rule{0cm}{33pt}\right.$}
\put(58,0){$\underbrace{\hskip19pt}_{\displaystyle c}$}
\put(228,0){$\underbrace{\hskip19pt}_{\displaystyle c}$}
\put(37,36){$d\left\{\rule{0cm}{33pt}\right.$}
\put(210,32){$b\left\{\rule{0cm}{19pt}\right.$}
\put(126,27){$*$}
\put(120,23){$\longrightarrow$}
\put(70,-28){{\bf (c)} Case $a> c$.}
\end{picture}
\bigskip\bigskip
     \includegraphics[scale=.6]{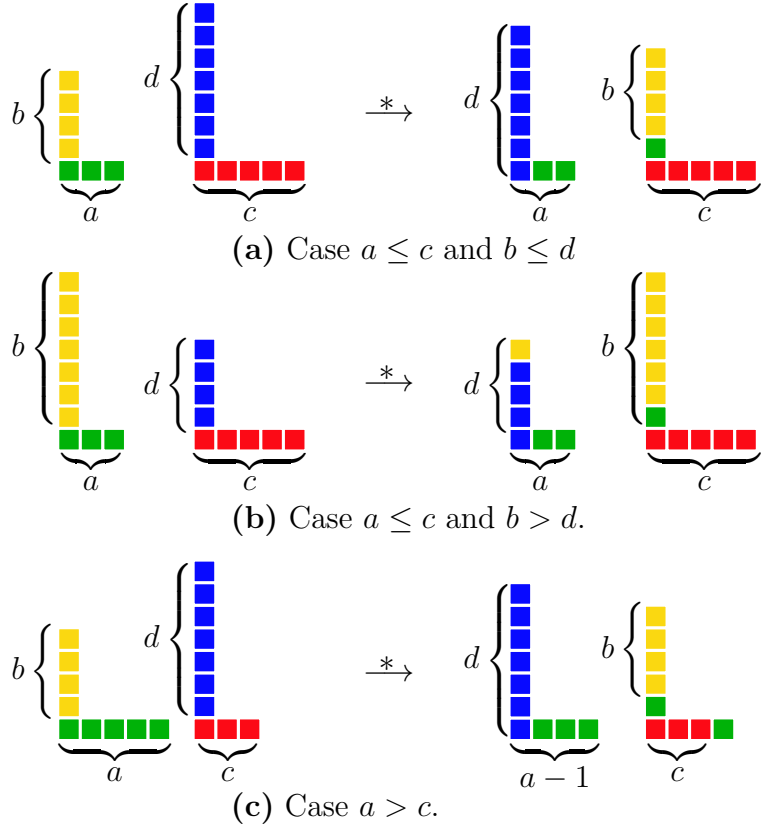}
    \caption{Three cases of $*$-operation on pairs of  hooks.}
    \label{star_hooks}
  \end{center}
\end{figure}

\goodbreak
We are now ready to prove Theorem \ref{theohooks} about other special instances of Conjecture \ref{fflp}, and of our extension Conjecture \ref{Bergeron}.
 \begin{proof}[\bf Proof of part (1) of  Theorem \ref{theohooks}] \quad
Let $\mu=(a,1^b)$ and $\nu=(c,1^d)$ be two hook shapes, and $(\lambda,\rho)$ be equal to $(\mu,\nu)^*$. There are essentially 3 different cases for the effect of  the $*$-operation on such a pair of hooks, depending on the relative values of $a$, $b$, $c$ and $d$. These are illustrated in Figure \ref{star_hooks}.  In each case, for $s_\theta$ that appears both in the expansion of $s_\mu s_\nu$ and $s_\lambda s_\rho$, our objective is to construct an injection between LR-fillings of $\theta/\nu$ of type $\mu$ and LR-fillings of $\theta/\rho$ of type $\lambda$. Under the above hypothesis, it is easy to check that $s_\theta$ can appear in the product $s_\mu s_\nu$, with nonzero coefficient, only if $\theta$ has at most two parts larger than $2$. Thus, in general, $\theta$ has the form $\theta=(r,s,2^t,1^u)$. Moreover, it is also clear that $r\geq \max(a,c)$, and $t+u+2\geq \max(b+1,d+1)$, since otherwise it would be impossible to get a nonzero result using the Littlewood-Richardson Rule.
\begin{figure}[ht] 
  \begin{center}
 \begin{picture}(0,0)(0,0)
 \put(-5,70){$d$}
  \put(185,57){$b$}
  \put(230,0){$c$}
  \put(45,0){$c$}
  \end{picture}
     \includegraphics[scale=.6]{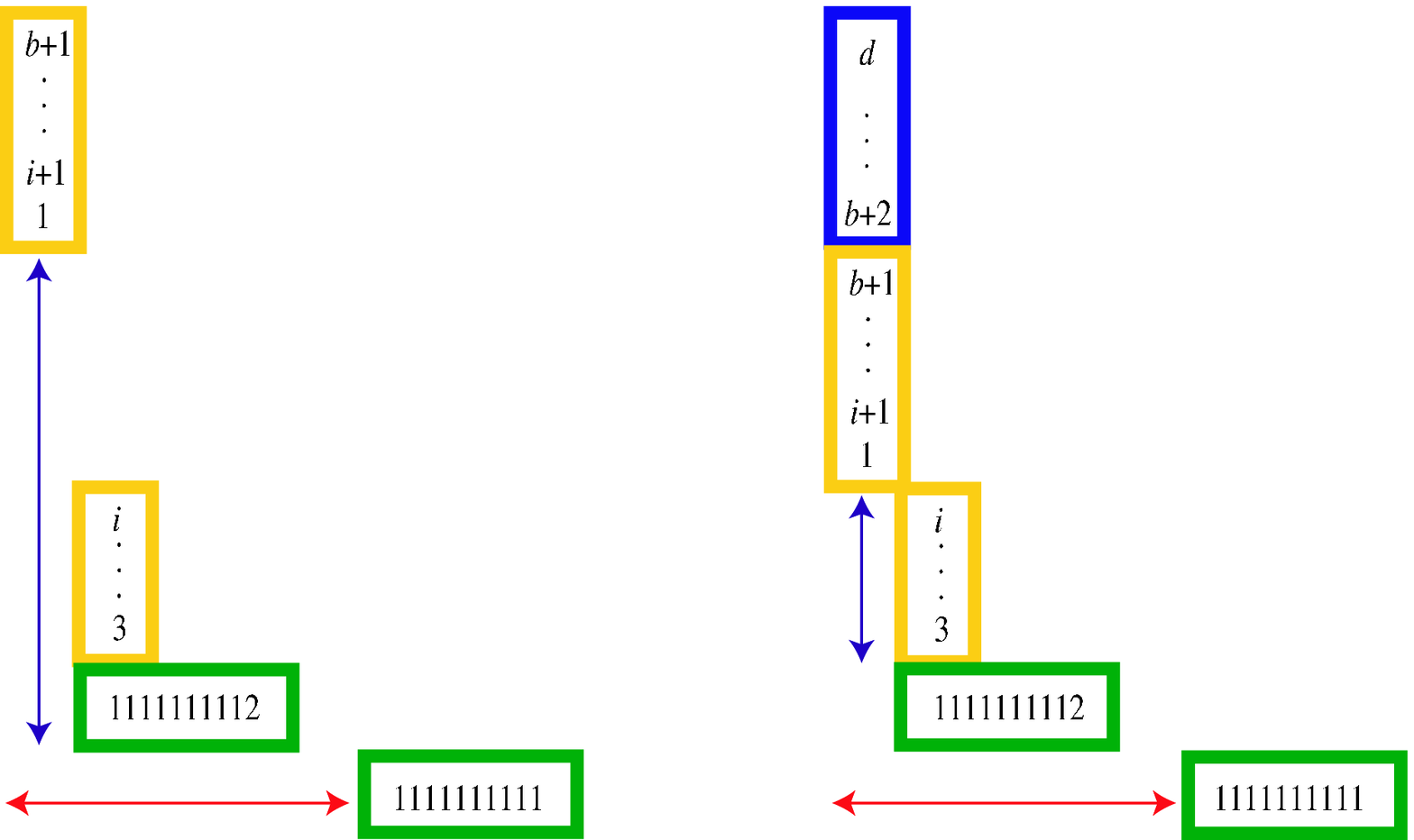}
    \caption{From a LR-fillings of $\theta/\nu$ to a LR-fillings of $\theta/\rho$}
    \label{casea}
  \end{center}
\end{figure}

  \item [\bf (a)]  ($a\leq c$ and $b\leq d$.) \quad
If $b=d$ or $b+1=d$, then $(\mu,\nu)$ is a fixed point and the result is obvious. We can thus suppose that $d\geq b+2$. This situation is illustrated in Figure \ref{star_hooks} case {\bf (a)}, and we have $\lambda=(a,1^{d-1})$ and $\rho=(c,1^{b+1})$.
Thus, the skew shape $\theta / \rho$ only  differs from that of $\theta / \nu$ in the first column. There are now $d-(b+1)$ new boxes to be filled, which are all at the top end of the first column of $\theta$. From a LR-filling of $\theta / \nu$, of type $\mu$, we construct a filling of $\theta / \rho$ as follows. We simply slide down by $d-(b+1)$ positions the entries appearing in the first column, and then add new entries $b+2,\ldots, d$ in the $d-(b+1)$ cells at the top of this first column. All the other entries of the original filling are kept as they were. The resulting LR-filling is clearly of type $\rho$.
An example of this procedure is given in Figure \ref{casea}.
 
 \goodbreak  \item [\bf (b)] ($a\leq c$ and $b> d$.) \quad
In a sense, this is almost the reverse of case {\bf (a)}, with a shift of values by 1. Thus, we can just reverse the process described above.

\goodbreak \item [{\bf (c)}] ($a\leq c$.)\quad 
In this case $\lambda=(a-1,1^{d-1})$ and  $\rho=(c+1,1^{b+1})$, so that this is almost as in case {\bf (a)} but with a small variation of $\pm 1$ in the first rows. We can thus proceed just as in case {\bf (a)} or {\bf (b)}, depending on the relative values of $b$ and $d$, removing a $1$ from the first row to get an LR-filling of the right shape.
\end{proof}

\begin{proof}[\bf Proof of part (2) of Theorem \ref{theo2rows}]\quad
It suffices to show that the conjecture holds for pairs of
two-row shapes, since Proposition \ref{cortranspose}
will then imply that the conjecture will also hold for pairs of two-column shapes.
Let $\mu=(a,b)$ and $\nu=(c,d)$ be a pair
of two row shapes. Assume  that $d>1$. If $d=1$
the calculations are very similar and are left to
the reader. There are six cases to consider corresponding to the six different  linear extensions of  the 
partial order $a \ge b,$ and $ c \ge d$.
 For instance, If $c \ge d > a \ge b$ then 
   $$( (a,b), (c,d) ) ^*= ( (a+1, b), (c, d-1) ).$$
 Given a shape $\theta$, together with a  LR-filling of 
$\theta/(c, d)$ of type $(a,b)$,
we construct a LR-filling of $\theta/(c, d-1)$
of type $(a+1, b)$ by placing a $1$ in the extra square. 
(See Figure \ref{tworow}.)
Note that since $d-1 \ge a$ the result is semistandard, and
by construction it is a reverse lattice permutation.
\begin{figure*}[ht]
\centering 
\begin{picture}(10,20)(0,0)
\put(22,39){$\scriptstyle (a,b)$}
\put(105,39){$\scriptstyle (c,d)$}
\put(102,-10){$\scriptstyle (c,d-1)$}
\put(20,-10){$\scriptstyle (a+1,b)$}
\put(265,7){$\theta$}
\end{picture}
\includegraphics[width=110mm]{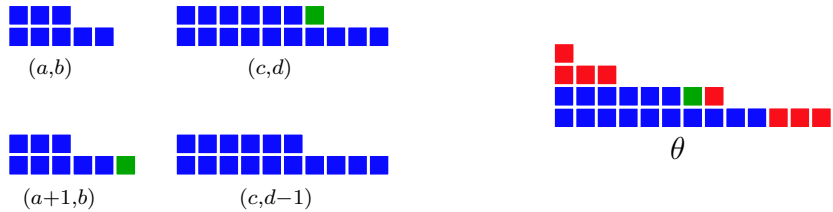} 
\bigskip
\caption{Case $c \ge d > a \ge b$.
\label{tworow}  }
\end{figure*}
 The other five cases  are all very similar. \end{proof}

\begin{proof}[\bf Proof of part  (3) of Theorem \ref{theo2rows}.]\quad To show that Conjecture \ref{Bergeron} holds for 
     $$(\mu/\alpha,\nu/\beta),$$
 when all partitions involved are hooks, we proceed as follows.
The typical case is when both $\alpha$ and $\beta$ are equal to the partition with one part of size $1$. It is easy to check that other situations are either similar or equivalent. Furthermore, as was done in  the proof of part {\bf (1)}, we will restrict our argument to the situation when $\mu=(a+1,1^b)$ and $\nu=(c+1,1^d)$, with $a\leq c$ and $d> b+1$. In that case, we have already observed that the result  $(\lambda,\rho)$ of the $*$-operation on $(\mu,\nu)$ is such that $\lambda=(a+1,1^{d-1})$ and $\rho(c+1,1^{b+1})$. It follows that, $(1,1)$ being a fixed point,  we have an explicit expression for $(\mu/1,\nu/1)^*$.
Furthermore (See \cite{macdonald}) we have
     $$s_{\mu/1}=h_a\,e_b,$$
 when $\mu=(a+1,1^b)$.
Here (as usual) $e_b$ denotes the classical {\em elementary} symmetric function. It follows that the identity that we have to check
is
\begin{eqnarray*}
    s_{\lambda/1} s_{\rho/1}-s_{\mu/1} s_{\nu/1} &=&
           h_a e_{d-1} h_c e_{b+1} - h_a e_b h_c e_d\\
           &=& h_a h_c (e_{d-1} e_{b+1} -e_d e_b).
\end{eqnarray*}
To conclude, we need only recall the ``dual'' Jacobi-Trudi identity and Pieri's formula to see that  last expression is Schur positive. The involution $\omega$ (introduced before Proposition \ref{cortranspose}) sends $h_a$ to $e_a$. Thus (the special case of) the Jacobi-Trudi identity
   $$s_{d-1,b+1}=\det \left( \begin{array}{cc}  
                                        h_{d-1} & h_d \\ 
                                        h_b     & h_{b+1} \end{array} \right)$$
becomes, after application of $\omega$,
   $$s_{2^{b+1},1^{d-b-2}}=\det \left( \begin{array}{cc}  
                                        e_{d-1} & e_d \\ 
                                        e_b     & e_{b+1} \end{array} \right).$$
To finish the proof, we recall that Pieri's formula underlines that $h_a s_\mu$ is Schur positive for any $a$ and any $\mu$.
\end{proof}

\begin{proof}[\bf Proof of part  (4) of Theorem \ref{theo2rows}.]\quad For $\nu/\beta$ a ribbon, set  $(\lambda,\rho)=(0,\nu/\beta)^*$. We will show that $s_\lambda s_\rho$ can be expanded as a (positive integer coefficient) sum of skew Schur functions indexed by ribbons, with $s_{\nu/\beta}$ appearing with nonzero coefficient. The slightly more general case of weak ribbons is entirely similar. We first need to recall that a $N$-cell ribbon $R=\nu/\beta$ is entirely described by the list of (positive) values $c_{j+1}:=\nu_{n-j}-\beta_{n-j}$, $0\leq j\leq n-1$, with $N=c_1+c_2+\ldots +c_n$. In other words, the $c_i$'s are the consecutive lengths of horizontal strips in $R$, starting from the top.
This list can itself be bijectively encoded as the subset of partial sums
     $$U(R):=\{c_1,c_1+c_2,c_1+c_2+c_3,\ldots, c_1+c_2+\ldots + c_{n-1} \},$$
of the set $\{1,\ldots, N-1\}$. Now, a classical result (See \cite{stanley}) about ribbon Schur functions states that
\begin{equation}\label{ribbon_expansion}
     \prod_{a\in A} h_a = \sum_{U(T)\subseteq U(A)} s_T,
\end{equation}
 where $A$ is any list of positive integers adding up to $N$, and $T$ varies among $N$-cells ribbons. This follows easily from the multiplication rule of ribbon Schur functions.
One can derive a formula analogous to Formula (\ref{ribbon_expansion}) which involves the elementary functions $e_e$ instead of $h_b$, with indices on the right hand side involving the length of vertical strips in $R$. A mixture of the two will be used below.
 
It follows readily from the definitions that the skew shapes $\lambda=\overline{\nu}/\overline{\beta}$ and
$\rho=\underline{\nu}/\underline{\beta}$ are such that
   $$s_\lambda = \prod_{a\in A} e_a \qquad {\rm and} \qquad s_\rho=  \prod_{b\in B} h_b,$$
where $A$ and $B$ are lists of the form
     $$A=(\underbrace{1,1,\ldots,1}_j,\ldots,d_2,d_1)\qquad {\rm and}\qquad B=(c_{i}-j,c_{i+1},\ldots,c_n),$$
 with the $d_i$'s consecutive lengths of vertical strips in $R$, starting from the top. More precisely, $A$ is the list of vertical lengths in the portion of the ribbon that lies in the region of $\nu$ that contributes to $\overline{\nu}$.
 This makes it evident that $s_R$ appears in the ribbon expansion of the product of $s_\lambda$ and $s_\rho$. Thus the theorem is proved.   
\end{proof}

\section{Reduction to a finite set of pairs in bounded height case}\label{reduction}

In this section we show that  the  bounded height case of Conjecture \ref{fflp} can be reduced to checking that it holds for a finite number of pairs, for any given height. In order to do this, and to state our result, we need some definitions.
Let $\mu$ and $\theta$ be two partitions such that $\mu \subseteq \theta$ and consider the skew partition $\theta / \mu$. Given a partition $\mu $ containing a column of height $k$,  we denote $\mu - 1^k$ the partition obtained by removing this column. In other words,
   $$(\mu -  1^k)_j :=\left \{\begin{array}{ll}
    \mu_j - 1, & \mbox{if} \; j \le k, \\ 
     \mu_j, & \mbox{if} \; j>k.
    \end{array} \right. $$
\begin{figure*}[ht]
\centering 
\includegraphics[ width=50mm]{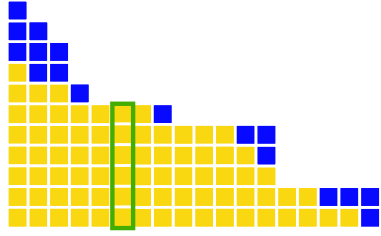} 
\caption{A 6-full column.}
\label{full} 
\end{figure*}
We say that $\mu$ has a  $k$-{\em full} column in $\theta$, if there is a $j$ such that  the $j^{\rm th}$ column of $\mu$ and  $\theta$ are both of height $k$.  When this is the case, setting $\beta:=\theta-1^k$ and $\gamma:=\mu-1^k$, we observe that
\begin{equation}\label{equal}
    \langle s_{\theta},s_{\mu}s_{\nu}\rangle=\langle s_{\beta},s_{\gamma }s_{\nu}\rangle,
\end{equation}
since, using (\ref{scalar}), this is clearly equivalent to $s_{\theta/\mu}=s_{\beta/\gamma}$ which holds trivially. When $\mu\subseteq \theta$, the fact that $\mu$ has a $k$-full column in $\theta$ is equivalent to
 \begin{equation}\label{criteria}
     \theta_k\geq  \mu_k>\theta_{k+1},
 \end{equation}
 assuming that $\theta_{k+1}=0$ when $k$ is the number of parts of $\theta$.

\begin{proof}[\bf Proof of  Theorem \ref{asympt}] \quad
\quad We proceed by induction on the number of columns of $\mu$, for the set of partitions $\mu$ with height bounded by $p$, and such that  
   \begin{equation}\label{assez}
       \mu_1 > p\,(\nu_1+p).
    \end{equation}
Once again, let $(\lambda,\rho):=(\mu,\nu)^*$. For any $\theta$, such that $s_\theta$ appears with nonzero coefficient in $s_\mu s_\nu$, we will show (\ref{assez}) implies  that there exists a $k$, such that both $\mu$ and $\lambda$ have a $k$-full column in $\theta$, and such that
\begin{equation}\label{stabilite}
    \lambda(\mu-1^k,\nu)=\lambda-1^k \qquad {\rm and}\qquad  \rho(\mu-1^k,\nu)=\rho.
  \end{equation}
It will then follow that
\begin{align*}
\langle s_{\theta},s_{\mu}s_{\nu} \rangle & =  \langle s_{\theta - 1^k},s_{\mu - 1^k}s_{\nu}\rangle && \text{by (\ref{equal})}\\
& \leq  \langle s_{\theta - 1^k},s_{\lambda(\mu - 1^k,\nu)}s_{\rho(\mu - 1^k,\nu)}\rangle && \text{by  induction hypothesis}\\
 & =  \langle s_{\theta - 1^k},s_{\lambda - 1^k }s_{\rho}\rangle  && \text{by (\ref{stabilite})}\\
 & =  \langle s_{\theta},s_{\lambda}s_{\rho}\rangle
 && \text{by (\ref{equal})}
\end{align*}
which will prove the theorem.

To show that there is a $k$ with the properties announced above, we proceed as follows. Observe that at least one of the differences $\mu_j-\mu_{j+1}$, where $1\leq j\leq p$, is strictly larger then $\nu_1+p$, since otherwise
\begin{eqnarray*}
    \mu_1&=&\sum_{j=1}^p \mu_j-\mu_{j+1}\\
               &\leq & p\,(\nu_1+p)
\end{eqnarray*}
which would contradict (\ref{assez}). We can thus choose $k$ to be the smallest integer, between $1$ and $p$, such that 
   \begin{equation}\label{choix_de_k}
          \mu_k>\mu_{k+1}+\nu_1+p.
  \end{equation}
For $\theta$ as above, we must clearly have $\mu_i\leq \theta_i\leq \mu_i+\nu_1$. Thus 
 \begin{eqnarray*}
       \theta_k&\geq & \mu_k\\
                     &> & \mu_{k+1}+\nu_1+p\\
                     &>& \theta_{k+1},
\end{eqnarray*}
so that $\mu$ has a $k$-full column in $\theta$, by criteria (\ref{criteria}). Moreover, 
for $1\leq i\leq k$, it is clear that $\mu_i>\nu_1+p$, and thus (\ref{definition}) simplifies to
    $$\lambda_i=\mu_i-i,\qquad {\rm for}\qquad 1\leq i\leq k.$$
 It follows that
\begin{eqnarray*}
    \lambda_k&=&\mu_k-k\\
                       &>&\mu_{k+1}+\nu_1+(p-k)\\
                       &\geq&\theta_{k+1},
 \end{eqnarray*}
 so that $\lambda$ also has a $k$-full column in $\theta$. The last verification that we need to do is that (\ref{stabilite}) holds. Now, it is clear that the first $k$ lines of $\gamma:=\mu-1^k$ are all too large for the first part of (\ref{gauche}) to apply. In fact, considering the way $k$ has been chosen, we see that for $1\leq i\leq k$
\begin{eqnarray*}
      \gamma_i-i&=&\mu_i-(i+1)\\
                  &>&\mu_1-(i-1)\, (\nu_1+p)-(i+1)\\
                  &>& (p-i+1)\,(\nu_1+p)-(i+1)\\
                  &\geq&\nu_1-1.
\end{eqnarray*}
This makes it obvious that (\ref{stabilite}) holds, thus finishing our proof.
\end{proof}

\section{Final remarks}\label{finalremarks}  

 We believe that to get a better understanding of the $*$-operation, a refined  study of its effect on tableaux and semistandard tableaux will be crucial. For instance this should lead to a proof of  ``monomial'' versions of Conjectures \ref{fflp} and \ref{Bergeron}. More precisely, recall that the expansion of any Schur function in the basis of monomial symmetric functions involves only positive integers. It would thus follow from the conjectures that the expansion of the difference of products considered have positive integer coefficients when expanded in term of monomial symmetric function. In particular, using definition (\ref{nombre_tableaux}), one should have 
  \begin{equation}
        f_{(\lambda,\rho)} \geq f_{(\mu,\nu)}.
  \end{equation}
  whenever $(\lambda,\rho)={(\mu,\nu)}^*$.
An independent proof of these facts would clearly lend support to the conjectures. 

\section{Acknowledgments}
We would like to thank Peter McNamara for the proof of Lemma \ref{dominance} presented in this paper, which is much nicer than our original proof. We are also grateful to Sergey Fomin for reading an earlier version of this manuscript, and for many useful suggestions.

\end{document}